%


\hfuzz=2pt
\def\sqr#1#2{{\vcenter{\vbox{\hrule height.#2pt \hbox{\vrule width.#2pt height#1pt \kern#1pt \vrule width.#2pt} \hrule height.#2pt}}}}

\def\qed{\hfill{\vbox{\hrule\hbox{\vrule\kern3pt
                \vbox{\kern6pt}\kern3pt\vrule}\hrule}}}
\def\Dom{{\rm Dom}}
\def\Ran{{\rm Ran}}
\def\cf{{\rm cf}}
\def\id{{\rm id}}

\def\ll{{\lambda}}
 
\def\chr{{\rm Chr}}
\def\ff{{\varphi}}
\def\dd{{\delta}}
\def\aa{{\alpha}}
\def\bb{{\beta}}

\def\tp{{\rm tp}}
\def\xx{\overline x_0,\dots,\overline x_n}
\def\yy{\overline y_0,\dots,\overline y_n}

\def\today{
   \ifcase\month\or January\or February\or March\or April\or May\or
    June\or July\or August\or September\or October\or
    November\or December\fi\space\number\day, \number\year}
\def\P{\noindent{\bf Proof.} }
\magnification1200

\centerline{\bf On Taylor's problem}
\medskip
\centerline{{\bf P.~Komj\'ath}\footnote{}{The first
author acknowledges the support of the Hungarian OTKA grant 2117.}{\bf\ and 
S.~Shelah}\footnote{}{Publication 
number No.~346 on the second author's list. His research was
partially covered by the Israel Academy Basic Research Fund.}}
\medskip
\def\KK{{\cal K}^{n,e}}
By extending finite theorems Erd\H os and Rado proved that for
every infinite cardinal $\kappa$ there is a $\kappa$-chromatic
triangle-free graph [3]. In later work they were able to add the
condition that  the graph itself be of cardinal $\kappa$ [4]. The next stage,
eliminating 4-circuits, turned out to be different, as it was
shown by Erd\H os and Hajnal [1] that every uncountably chromatic
graph contains a 4-circuit. In fact, every finite bipartite
graph must be contained, but odd circuits can be omitted up to a
certain length. This solved the problem ``which finite graphs
must be contained in every $\kappa$-chromatic graph'' for every
$\kappa> \omega$. The next result was given by Erd\H os, Hajnal,
and Shelah [2], namely, every uncountably chromatic graph
contains all odd circuits from some length onward. They, as well
as Taylor, asked the following problem. If $\kappa$, $\lambda$
are uncountable cardinals and $X$ is a $\kappa$-chromatic graph,
 is there a $\lambda$-chromatic graph $Y$ such that every finite
subgraph of $Y$ appears as a subgraph of $X$. In [2] the
following much stronger conjecture was posed. If $X$ is
uncountably chromatic, then for some  $n$ it contains all finite
subgraphs of the so-called $n$-shift graph. This conjecture was,
however, disproved in [5]. 

Here we give some results on Taylor's conjecture when the
additional hypotheses $\vert X\vert =\kappa$, $\vert Y\vert
=\ll$ are imposed. 

We describe some (countably many) classes  $\KK$ of finite
graphs and prove that if $\ll^{\aleph_0}=\ll$ then every
$\ll^+$-chromatic graph of cardinal $\ll^+$ contains, 
for some $n$, $e$, all members
of $\KK$ as subgraphs. On the other hand, it is consistent for
every regular infinite cardinal $\kappa$ that there is a
$\kappa^+$-chromatic graph on $\kappa^+$ that contains finite
subgraphs only from $\KK$. We get, therefore, some models of set
theory, where the finite subraphs of graphs with $\vert X\vert=
\chr(X) =\kappa^+$ for regular uncountable cardinals $\kappa$
are described. 

We notice that in [6] all countable graphs are described which
appear in every graph with uncountable coloring number.

\medskip
\noindent{\bf Notation.} $\overline x$ will denote  a finite string of
ordinals. $\overline x < \overline y$ means that $\max
(\overline x) < \min (\overline y)$.
\medskip
\noindent{\bf Definition.} Assume that $1\leq n<\omega $, $e:
\{1,2,\dots,2n\}\to \{0,1\}$ is a function with $\vert
f^{-1}(0)\vert =n$. We are going to define the structures in
$\KK$ as follows. They will be of the form
$H=(V,<,U,X,h_1,\dots,h_n)$ where $(V,<)$ is a finite linearly ordered
set, $U\subseteq V$, $X$ is a graph on $U$, $h_i:U\to V$ satisfy
$h_1(x)<\cdots<h_n(x)=x$ for $x\in U$. The elements in $\KK_0$
are those isomorphic to $(V,<,U,X,h_1,\dots,h_n)$ where
$V=\{1,2,\dots,n\}$, $<$ is the natural ordering, $U=\{n\}$,
$X=\emptyset$, $h_i(n)=i$ ($1\leq i\leq n$).

If $H=(V,<,U,X,h_1,\dots,h_n)$ is a structure of the above form,
and $x\in V$, we form the edgeless amalgamation $H'=H+_x H$ as
follows.  Put $H'=H +_x H=(V',<',U',X',h'_1,\dots,h'_n)$ where
$(V',<')$ has the {\sl $<'$-ordered} decomposition $V'=W\cup
V_0\cup V_1$, if we put $V'_i=W\cup V_i$ for $i<2$ then the
structures 
$$
\bigl(V'_i,<'\vert V'_i,U'\cap V'_i,h'_1\vert V'_i,\dots,h'_n\vert V'_i\bigr)
$$ 
are both isomorphic to $H$ for $i<2$ and $\min(V_i)$ correspond to
$x$ under the isomorphisms.

If $H=(V,<,U,X,h_1,\dots,h_n)$ is a structure of the above form,
and $x\in U$, we also form the one-edge amalgamation $H'=H*_x H$ as
follows. Enumerate in increasing order $e^{-1}(0)$ as
$\{a_1,\dots,a_n\}$ and $e^{-1}(1)$ as
$\{b_1,\dots,b_n\}$. Put $H'=(V',<',U',X',h'_1,\dots,h'_n)$ where
$(V',<')$ has the {\sl ordered} decomposition $V'=V_0\cup
V_1\cup\cdots\cup V_{2n}$; $H'\vert (V_0\cup\bigcup\{V_i:e(i)=
\varepsilon\})$ are isomorphic to $H$ ($\varepsilon=0,1$) if
$x_0$, $x_1$ are the points 
corresponding to $x$, then $h'_i(x_0)=\min(V_{a_i})$, 
$h'_i(x_1)=\min(V_{b_i})$, and the only extra edge in $X'$ is $\{x_0,x_1\}$.

We then put 
$$
\eqalign{\KK_{t+1}=&\bigl\{H+_x H:H=(V,<,U,\dots)\in \KK_t, x\in
V\bigr\}
\cr&\cup  \bigl\{H*_y H:H=(V,<,U,\dots)\in \KK_t, y\in
U\bigr\},\cr}
$$
and finally $\KK=\bigcup\{\KK_t:t<\omega \}$.

\medskip
\proclaim Theorem 1. If $\vert G\vert =\chr(G)=\ll^+$,
$\ll^{\aleph_0}=\ll$, then, for some $n$, $e$, $G$ contains 
every graph in $\KK$ as subgraph.

\medskip
We start with some technical observations.
\medskip
\proclaim Lemma 1. If $t_n:\ll^+\to\ll^+$ are functions
($n<\omega $), then there is a $\ll$-coloring $F:\ll^+\to\ll$
such that for $F(\alpha )=F(\beta )$, $i$, $j<\omega $, $\alpha
<t_i(\beta )<t_j(\alpha )$ may not hold.

\medskip
\P As $\ll^{\aleph_0}=\ll$, it suffices to show this for two
functions $t_0(\alpha )$, $t_1(\alpha )$, with $t_1(\alpha
)>\alpha $. We prove the stronger statement that there is  a function
$F:\ll^+\to[\ll]^\ll$ such that if $\alpha <t_0(\beta
)<t_1(\alpha )$ then $F(\alpha )\cap F(\beta )=\emptyset$. Let 
$\langle N_\xi:\xi<\ll^+\rangle$ be a continuous, increasing
sequence of elementary submodels of
$\langle\ll^+;<,t_0,t_1,\dots\rangle$ with $\gamma
_\xi=N_\xi\cap \ll^+<\ll^+$. $C=\{\gamma _\xi:\xi<\ll^+\}$ is
closed, unbounded. We define $F\vert \gamma _\xi$ by transfinite
recursion on $\xi$. If $F\vert \gamma _\xi$ is given, and $\beta
$ has $t_0(\beta )<\gamma _\xi\leq \beta <\gamma _{\xi+1}$, by
elementarity $\tau=\sup\{t_1(\alpha ):\alpha <t_0(\beta
)\}<\gamma _\xi$ , and there is a $\beta '$ with $t_0(\beta
')=t_0(\beta )$,  $\tau<\beta '<\gamma _\xi$. Put $H(\beta
)=F(\beta ')$, otherwise, i.e., when $\gamma_\xi\leq t_0(\beta
)$, put $H(\beta )=\ll$. To get $F\vert [\gamma _\xi,\gamma
_{\xi+1})$, we disjointize $\{H(\beta ):\gamma _\xi\leq \beta
<\gamma _{\xi+1}\}$, i.e., find $F(\beta )\subseteq H(\beta )$
of cardinal $\ll$ 
such that $F(\beta _0)\cap F(\beta _1)=\emptyset$ for $\beta
_0\neq \beta _1$. 
We show that  this $F$ works. 
Assume  that $F(\alpha )$, $F(\beta )$ are not disjoint. By
induction we can assume that either $\alpha$ or $\beta$ is
between $\gamma_\xi$ and $\gamma_{\xi+1}$. By the disjointization
process some of them must be smaller than $\gamma_\xi$. If
$\beta<\gamma_\xi\leq \alpha<\gamma_{\xi+1}$ then
$t_0(\beta)<\gamma_\xi$ as $N_\xi$ is an elementary submodel,
so $t_0(\beta)<\alpha$. Assume now that 
$\alpha
<\gamma _\xi\leq \beta <\gamma _{\xi+1}$.
Our construction then
selected a $\beta '$ with $t_0(\beta ')=t_0(\beta )$ and
$F(\beta )\subseteq H(\beta )=F(\beta ')$ which is, by the
inductive hypothesis, disjoint from
$F(\alpha )$. \qed

\medskip

\medskip
\noindent{\bf Lemma 2.} {\sl If $C=\{\delta _\xi:
\xi<\ll^+\}$ is a club then there is a function
$K:[\ll^+]^{\aleph_0} \to \ll$ such that if $K(A)=K(B)$ and 
$A\cap [\delta _\xi,\delta _{\xi+1})\neq\emptyset$ and 
$B\cap [\delta _\xi,\delta _{\xi+1})\neq\emptyset$ for some
$\xi<\ll^+$ 
then $A\cap \delta _{\xi+1}=B\cap \delta _{\xi+1}=$ and 
so $A\cap B$ is an initial segment  both in $A$ and $B$.}
\medskip
\P 
Fix for every $\beta <\ll^+$ an into function 
$F_\beta :\alpha \to \ll$ such that for $\beta _0<\beta _1<\beta
_2$, $F_{\beta  _1}(\beta _0)\neq F_{\beta _2}(\beta _1)$ holds.
This can be done by a straightforward inductive construction. 

If $A\in [\ll^+]^{\aleph_0}$ put $X(A)=\{\xi:
A\cap [\delta _\xi,\delta _{\xi+1})\neq\emptyset\}$. 
Let $\tp(X(A))=\eta$. 
Enumerate $X(A)$ as $\{\tau^A_\theta:\theta<\eta\}$. 
Let $K(A)$ be a function with domain $\eta$, at $\theta<\eta$,
if $\tau^A_\theta=\xi$, let 
$$
K(A)(\theta)=\big\langle\{F_{\tau^A_\theta}(\tau^A_{\theta'}):
\theta'<\theta\},
\{F_{\delta _{\xi+1}}(y):y\in A\cap \delta _{\xi+1}\}\big\rangle.
$$
Assume now that $K(A)=K(B)$, $\xi\in X(A)\cap X(B)$. If $\xi=\tau^A_\theta=
\tau^B_{\theta'}$ then $\theta=\theta'$ by the properties of $F$
above. The second part of the definition of $K(A)$ gives that
$A\cap \delta _{\xi+1}=B\cap \delta _{\xi+1}$.
 \qed

\medskip
\noindent{\bf Proof of  Theorem 1.} 
We first show that one can assume that $G$ is $\ll^+$-chromatic
on every closed unbounded set.

\medskip
\noindent{\bf Lemma 3.} {\sl There is a function $f:\ll^+\to
\ll^+$ such that if $C\subseteq \ll^+$ is a closed unbounded set
then $\bigcup\bigl\{[\alpha ,f(\alpha )]:\alpha \in C\bigr\}$ 
is $\ll^+$-chromatic.}
\medskip
\P Assume that the statement of the Lemma fails. Put $f_0(\alpha
)=\alpha $, for $n<\omega$ let $C_n$ witness that
$f_n:\ll^+\to\ll^+$ is not good and $f_{n+1}(\alpha
)=\min(C_n-(\alpha +1))$. As, by assumption,  
$\bigcup\{[\alpha ,f_n(\alpha )]:\alpha \in C_n, n<\omega\}$ is
$\leq\ll$-chromatic, there is a 
$$
\gamma \notin 
\bigcup\Bigl\{\bigl[\alpha ,f_n(\alpha )\bigr]:
\alpha \in C_n,n<\omega \Bigr\},
\gamma>\min\Bigl(\bigcap\Bigl\{C_n:n<\omega\Bigr\}\Bigr).
$$ 
Clearly, $\gamma \notin C_n$
($n<\omega $), and if now $\alpha _n=\max(\gamma \cap C_n)$,
then $\alpha _n<\gamma $, and $\alpha _{n+1}<\alpha _n$
($n<\omega $), a contradiction. \qed

\medskip
By slightly re-ordering $\ll^+$ we can state Lemma 3 as follows.
If $C\subseteq \ll^+$ is a closed unbounded set, then $S(C)=
\bigcup \{[\ll\alpha ,\ll(\alpha +1)):\alpha \in C\}$ is
$\ll^+$-chromatic. Put, for $\tau<\ll$, $C\subseteq \ll^+$ a
club set, 
 $S_\tau(C)=\bigcup\{\ll\alpha
+\tau:\alpha \in C\}$. If, for every $\tau<\ll$ there is some
closed unbounded $C_\tau$ that $S_\tau(C_\tau)$ is
$\ll$-chromatic, then for $C=\bigcap\{ C_\tau:\tau<\ll\}$, $S(C)$
is the union of at most $\ll$ graphs, each $\leq\ll$-chromatic,
a contradiction. 

There is, therefore, a $\tau<\ll$ such that $S_\tau(C)$ is
 $\ll^+$-chromatic whenever $C$ is a closed unbounded set.
Mapping $\ll\alpha +\tau$ to $\alpha $ we get a graph on
$\ll^+$, order-isomorphic to a subgraph of the original graph
which is $\ll^+$-chromatic on every closed unbounded set. From
now on we assume that our original graph $G$ has this property.
\medskip
 We are going
to build a model $M=\langle\ll^+;<,\ll,G,\dots\rangle$ by adding
countably many new functions. 

For $n$, $e$ as in the Definition, $\ff$ a first order formula,
let $G^{n,e}_\ff$ be the following graph. The vertex set is
$V_\ff= \{\langle\overline x_0,\dots,\overline x_n\rangle:
\overline x_0<\cdots<\overline x_n,M\models\ff(\overline
x_0,\dots,\overline x_n)\}$ and $\langle\overline
x_0,\dots,\overline x_n\rangle$, $\langle\overline y_0,\dots,\overline 
y_n\rangle$ are joined, if $\overline x_0=\overline y_0$,
$\{\overline x_1,\dots,\overline x_n\}$,
$\{\overline y_1,\dots,\overline y_n\}$ interlace by $e$, and
finally $\{\min(\overline x_n),\min (\overline y_n)\}\in G$. We
introduce a new quantifier $Q^{n,e}$ with $Q^{n,e}\ff$ meaning
that the above graph, $G^{n,e}_\ff$ is $\ll^+$-chromatic. If,
however, $\chr(G^{n,e}_\ff)\leq\ll$, we add a good
$\ll$-coloring to $M$. We also assume that $M$ is endowed with
Skolem functions. 
\medskip
\proclaim Lemma 4. There exist $n$, $e$ and $\alpha
_1<\cdots<\alpha _n<\ll^+$ such that 
$t(\alpha _i)<\alpha _{i+1}$ if $t:\ll^+\to \ll^+$ is a function
in $M$ and if $\overline x_0\subseteq
\alpha _1$, $\overline x_i\subseteq [\alpha _i,\alpha _{i+1})$
($1\leq i< n$), $\overline x_n\subseteq [\alpha _n,\ll^+)$,
$\min(\overline x_i)=\alpha _i$, and $\ff$ is a formula,
$M\models \ff(\xx)$, then $M\models Q^{n,e}\ff$.

\medskip
\P Assume that the statement of the Lemma does not hold, i.e.,
for every $n$, $e$, $\alpha _1,\dots,\alpha _n$ there exist
$\xx$ contradicting it. 

Let, for $\alpha <\ll^+$,  $B_\alpha \subseteq \ll^+$ be a
countable set such that $\alpha \in B_\alpha $, and if $n$, $e$,
$\alpha _1,\dots,\alpha _n\in B_\alpha $ are given , then a
counter-example as above is found with $\xx\subseteq B_\alpha $.
We require that $B_\alpha $ be Skolem-closed. Let $B^+_\alpha $
be the ordinal closure of $B_\alpha $, $B^+_\alpha =\{\gamma
(\alpha ,\xi):\xi\leq \xi_\alpha \}$ be the increasing
enumeration, $\alpha =\gamma (\alpha ,\tau_\alpha )$. 
Let $\{M_\xi:\xi<\ll^+\}$ be a continuous, increasing chain of
elementary submodels of $M$ such that $\delta _\xi=M_\xi\cap
\ll^+<\ll^+$. Clearly, $C=\{\delta _\xi:\xi<\ll^+\}$ is a
closed, unbounded set. We take a coloring of the sets
$\{B^+_\alpha :\alpha <\ll^+\}$ by $\ll$
colors that satisfies Lemma 2, if $\alpha $, $\beta $ get the
same color then the structures $(B^+_\alpha ;B_\alpha ,M)$ 
and $(B^+_\beta ;B_\beta ,M)$ are isomorphic and we also
 require that if
$\xx\subseteq B_\alpha $ and $\yy\subseteq B_\beta $ are in the
same positions, i.e., are mapped onto each other by the order
isomorphism between $B_\alpha $ and $B_\beta $ and $(\xx)$ is
colored by the $\ll$-coloring of $G^{n,e}_\ff$, then $(\yy)$ is
also colored and gets the same color. All this is possible, as
$\ll^{\aleph_0} =\ll$. We also assume that our coloring
satisfies Lemma 1 with some functions 
$\{t_n:n<\omega\}$ that $B^+_\alpha=\{t_n(\alpha):n<\omega\}$.

\medskip
As $G$ is $\ll^+$-chromatic on $C$, there
are $\alpha <\beta $, both in $C$, joined in $G$, getting the same color.
%
By our conditions, $B^+_\alpha \cap B^+_\beta $ is initial
segment in both, and beyond that they do not even intersect into
the same complementary interval of $C$. As our structures are
isomorphic, this  holds for $B_\alpha $, $B_\beta $, as
well. 
\medskip
We now let $B^+_\alpha =\bigcup\{B^+_\alpha  (i):i<i_\alpha \}$, 
$B^+_\beta =\bigcup\{B^+_\beta  (i):i<i_\beta \}$ 
be the ordered decompositions given by the 
following equivalence relations. For $x$, $y\in
B^+_\alpha $, $x\leq y$, $x\sim y$ if either  $[x,y]\cap
B^+_\beta = \emptyset$  or  $[x,y]\cap B^+_\beta \supseteq 
[x,y]\cap B^+_\alpha$. 
Similarly for $B^+_\beta $. By
Lemma 2, $B^+_\alpha \cap B^+_\beta =B^+_\alpha (0)=B^+_\beta (0)$.

\medskip
\proclaim Lemma 5. $i_\alpha $, $i_\beta $ are finite.

\medskip
\P Otherwise, as $B^+_\alpha $, $B^+_\beta $ are ordinal closed,
$\gamma =\min\bigl(B^+_\alpha (\omega )\bigr)=
\min\bigl(B^+_\beta (\omega )\bigr)$ is in $B^+_\alpha \cap
B^+_\beta $, so $\gamma \in B^+_\alpha (0)$, a contradiction. \qed

\medskip
Enumerate 
$$\eqalign{
\Bigl\{\xi\leq\xi_\alpha :&\hbox{ there is a }0<i<\omega {\rm\
such\ 
that\ either\ }\gamma (\alpha ,\xi)=\min
\bigl(B^+_\alpha (i)\bigr)\cr &{\rm\ or\ }
\gamma (\beta  ,\xi)=
\min\bigl(B^+_\beta  (i)\bigr)\Bigr\}\cup\Bigl\{\tau_\alpha 
\Bigr\}\cr}
$$
as $\xi_1<\xi_2<\cdots<\xi_n$. By Lemma 1, if $\alpha <\beta $, 
$\alpha =\min(B^+_\alpha (i_\alpha -1))$, 
$\beta  =\min(B^+_\beta  (i_\beta  -1))$, $B^+_\alpha (i_\alpha
-1)< B^+_\beta (i_\beta -1)$. So $\xi_n=\tau_\alpha $. We let
$\alpha _i=\gamma (\alpha ,\xi_i)$, $\beta _i=\gamma (\beta
,\xi_i)$. 
If $\alpha _i=\min(B^+_\alpha (j))$ then $\alpha _i\in B_\alpha
$ and, by isomorphism, $\beta _i\in B_\beta $. 
We show that for every $i<n$, $t\in M$, $t(\alpha
_i)<\alpha 
_{i+1}$ and $t(\beta _i)<\beta _{i+1}$. As $(B^+_\alpha
;B_\alpha ,M)$ 
and $(B^+_\beta ;B_\beta ,M)$ are isomorphic, for every $i$ it
suffices to show this either for $\alpha _i$ or for $\beta _i$.
For $i=n-1$ this follows from the fact that $\alpha $ (as well
as $\beta $) is from $C$. If $i<n$ then either $\alpha _{i-1}$
and $\alpha _i$ are separated by an element of $B^+_\beta $ or
vice versa. Assume the former. Then, by Lemma 2, $\alpha _{i-1}$
and $\alpha _i$ are in different intervals of $C$ so necessarily
$t(\alpha _{i-1})<\alpha _i$ holds.

 Let $e$ be the interlacing type of $\{\alpha _i:1\le
i\le n\}$, $\{ \beta _i:1\le i\le n\}$. By our indirect
assumption, there are a formula $\ff$, $\overline x_i$, $\overline y_i$
($0\le i\le n$) in the same position in $B_\alpha $, $B_\beta $
such that $\overline x_i\subseteq [\alpha _i,\alpha _{i+1})$,
$\overline y_i\subseteq [\beta  _i,\beta  _{i+1})$ etc, and
$M\models \ff(\xx)\wedge\ff(\yy)$ and $(\xx)$, $(\yy) $ are joined
in $G^{n,e}_\ff$, but they get the same color in the good
coloring of $G^{n,e}_\ff$, a contradiction which proves Lemma 4. \qed
\medskip
Now fix $n$, $e$, and $\alpha _1<\cdots<\alpha _n<\ll^+$ as in
Lemma 4. We call a formula $\ff$ {\sl dense} if there exist 
$\overline x_0\subseteq
\alpha _1$, $\overline x_i\subseteq [\alpha _i,\alpha _{i+1})$
($1\leq i< n$), $\overline x_n\subseteq [\alpha _n,\ll^+)$,
$\min(\overline x_i)=\alpha _i$ such that $M\models \ff(\xx)$. 
If $H=(V,<,U,X,h_1,\dots,h_n)\in \KK$, $V=\{0,1,\dots,s\}$, a 
{\sl $\ff$-rich copy of $H$} is some string 
$(\overline y_0,\dots,\overline y_s)$ such that
$\overline y_0<\cdots<\overline y_s$,  if $\{i,j\}\in X$
then $\{\min(\overline y_i),\min(\overline y_j)\}\in G$ and for
every $v\in U$, $M\models \ff(\overline y_0,\overline y_{h_1(v)},
\dots,\overline
y_{h_n(v)})$.

\medskip
\proclaim Lemma 5. For every $H\in \KK$ if $\ff$ is dense there
is a $\ff$-rich copy of $H$ in $G$.

\medskip
\proclaim Lemma 6. For every $H=(V,<,U,X,h_1,\dots,h_n)\in
\KK$, $q\in U$, if $\ff$ is dense, there is a $\ff$-rich copy 
$(\overline y_0,\dots,\overline y_s)$ of $H$ such that
$\min(\overline y_{h_i(q)})=\alpha _i$ for $1\leq i\leq n$.

\medskip
We notice that Lemma 5 obviously concludes the proof of Theorem 1 and
Lemma 6 clearly implies Lemma 5. Also, they trivially hold for
$H\in \KK_0$. We prove these two Lemmas
simultaneously.

\medskip
\proclaim Claim 1. If Lemma 5 holds for some $H$ then Lemma 6
holds for $H$, as well.

\medskip
\P  Assume that Lemma 5 holds for 
$H=(V,<,U,X,h_1,\dots,h_n)\in
\KK$
 and for any dense $\ff$
but Lemma 6 fails for a certain $q\in U$ and a dense $\ff$. This
statement can be written as a formula 
$\theta(\alpha_1,\dots,\alpha _n)$. As $\ff $ is dense, 
$M\models \ff(\xx)$ for some appropriate strings, so also 
$M\models\psi(\xx)$ where $\psi=\ff\wedge 
\theta(\min(\overline x_1),\dots,\min(\overline x_n))$. As
$\psi$ is dense, by Lemma 5 there is a $\psi$-rich copy 
$(\overline y_0,\dots,\overline y_s)$ of $H$ but then 
$M\models \theta(\min(\overline y_{h_1(q)}),\dots,
\min(\overline y_{h_n(q)}))$, a contradiction. \qed

\medskip
\proclaim Claim 2. If Lemma 6 holds for
$H=(V,<,U,X,h_1,\dots,h_n)$ and $x\in V$ then Lemma 5 holds for
$H'=H +_x H$.

\medskip
\P Select $q\in U$ such that $x=h_i(q)$ for some $1\leq i\leq
n$. By Lemma 6, there is a $\ff$-rich copy 
$(\overline y_0,\dots,\overline y_s)$ of $H$ such that 
$\min(\overline y_{h_i(q)})=\alpha _i$ for $1\leq i\leq n$. 
As $\alpha _i>t(\alpha _{i-1})$ holds for every function $t$ in
the skolemized structure $M$ there
are $\ff$-rich copies of $H$ which agree with this below $x$ but
their $x$ elements are arbitrarily high. We can, therefore, get
a $\ff$-rich copy of $H'$. \qed

\medskip
\proclaim Claim 3. If Lemma 6 holds for
$H=(V,<,U,X,h_1,\dots,h_n)$ and $y\in U$ then Lemma 5 holds for
$H'=H *_y H$.

\medskip
\P Let $(\overline y_0,\dots,\overline y_s)$ be a $\ff$-rich
copy of $H$ such that $\min(\overline y_{h_i(q)})=\alpha _i$ 
for $1\leq i\leq n$. The elements in the $(\overline y_0,
\dots,\overline y_s)$ string can be redistributed as 
$(\overline x_0,\dots,\overline x_n)$ such that $\min(\overline
x_i)=\alpha _i$  and then the fact that they form a $\ff$-rich
copy of $H$ can be written as 
$M\models \psi(\xx)$ for some formula   $\psi$. As $\psi$ is
dense, by Lemma 4, $M\models Q^{n,e}\psi$ holds, so there are
two strings, 
$(\overline x_0,\dots,\overline x_n)$      and 
$(\overline x'_0,\dots,\overline x'_n)$ both satisfying $\psi$,
interlacing by $e$, and $\{\min(\overline x_n),\min(\overline
x'_n)\}\in G$. This, however, gives a $\ff$-rich copy of $H'$. \qed

\bigskip
\proclaim Theorem 2. If $n$, $e$ are as in the Definition, $\ll$
is an infinite cardinal, $\ll^{<\ll}=\ll$, then there exists a
$\ll^+$-c.c., $<\ll$-closed poset $Q=Q_{n,e,\ll}$ which adds a
$\ll^+$-chromatic graph of cardinal $\ll^+$ all whose finite
subgraphs are subgraphs of some element of $\KK$.

\P Put $q=(V,U,X,h_1,\dots,h_n)\in Q$ if $V\in
[\ll^+]^{<\ll}$, $U\subseteq V$, $X\subseteq [V]^2$, every $h_i$
is a function $U\to V$ with $h_1(x)<\cdots<h_n(x)=x$ for $x\in
U$ and every finite substructure of $(q,<)$ is a substructure of
some element of $\KK$. Order $Q$ as follows.
$q'=(V',U',X',h'_1,\dots,h'_n) \leq q=(V,U,X,h_1,\dots,h_n)$ iff
$V'\supseteq V$, $U=U'\cap V$, $X=X'\cap [V]^2$, $h'_i\supseteq
h_i$ ($1\leq i\leq n$). Clearly, $(Q,\leq)$ is $<\ll$-closed.

\medskip
\proclaim Lemma 7. $(Q,\leq)$ is $\ll^+$-c.c.

\P  By the usual $\Delta$-system arguments it suffices to show
that if the conditions 
$q^i=(V\cup V^i,U^i,X^i,h^i_1,\dots,h^i_n)$ are order
isomorphic ($i<2$), $V< V^0 < V^1$ then they are compatible. A
finite subset of $V\cup V^0 \cup V^1$ can be included into some 
$s\cup s_0 \cup s_1$ where $s_0$ and $s_1$ are mapped onto each
other by the isomorphism between $q_0$ and $q_1$. By condition, 
$q \vert s \cup s_0$ is a substructure of some structure $H\in
\KK$. But then $q\vert s\cup s_0 \cup s_1$ is a substructure of
an edgeless amalgamation of $H$. 
\qed

\medskip
If $G\subseteq Q$ is generic then $Y=\bigcup\{X:(V,U,X,\dots)\in
G\}$ is a graph on a subset of $\ll^+$ all whose finite
subgraphs  are subgraphs of some member of $\KK$. The
following Lemma clearly  concludes the proof of the Theorem.

\medskip
\proclaim Lemma 8. $\chr(Y)=\ll^+$.

\medskip
\P Assume, toward a contradiction, that 1 forces that
$f:\ll^+\to \ll$ is a good coloring of $Y$. Let $M_1\prec
M_2\prec \cdots \prec M_n$ be elementary submodels of
$(H((2^\ll)^+);Q,f,\dots ,)$ with $\ll \subseteq M_0$,
$[M_i]^{<\ll}\subseteq M_i$. Put $\dd_i=M_i\cap \ll^+< \ll^+$.
Notice that $\cf(\dd_i)=\ll$. Let $p'=(V',U',X',h'_1,\dots,h'_n)$
where $V=\{\dd_1,\dots,\dd_n\}$, $U=\{\dd_n\}$, $X=\emptyset$,
$h_i(\dd_n)= \dd_i$. Choose $p=(V,U,X,h_1,\dots,h_n)\leq p'$ 
forcing $f(\dd_n)=\xi$
for some $\xi<\ll$. Let $\psi_n(\pi,x_1,\dots,x_n)$ be the
following formula. $\pi$ is an order isomorphism $V\to \ll^+$,
$\pi(\dd_i)= x_i$ 
and $\pi(p)$ forces that $f(x_n)=\xi$. Let $\dd_{n+1}=\ll$. For
$0\leq i<n$ define $\psi_i(\pi,x_1,\dots,x_i)$ meaning that 
$\pi:V\cap \dd_{i+1}\to \ll^+$ is order preserving and 
there are arbitrarily large $x_{i+1}<\ll^+$ and $\pi'\supseteq \pi$
such that $\psi_{i+1}(\pi',x_1,\dots,x_{i+1})$ holds. 

\medskip \proclaim Claim 4. $\psi_i(\id\vert V\cap \dd_{i+1},
\dd_1,\dots, \dd_i)$ for $0\leq i\leq n$.

\medskip
\P This is obvious for $i=n$. If $\psi_i(\id\vert V\cap \dd_{i+1},
\dd_1,\dots, \dd_i)$ fails, then, by definition, there would be
a bound for the possible $x_{i+1}$ values for which 
$\psi_{i+1}(\pi',\dd_1,\dots, \dd_i,x_{i+1})$ holds for some
$\pi'\supseteq \id\vert V\cap \dd_{i+1}$. But then this bound is
smaller than $\dd_{i+1}$ so $\psi_{i+1}$ fails, too. 
\qed 

\medskip
Returning to the proof of Lemma 8, we define the following
function $t$. Let $\{a_1,\dots,a_n\}$, $\{b_1,\dots,b_n\}$ be as in
the definition of the one-edge amalgamation. Put, for $1\leq
i\leq n$, $t(i)=j$ iff $b_{j-1}<a_i<b_j$ where $b_0=0$,
$b_{n+1}=2n+1$. Set $\pi_0=\id\vert V\cap \dd_1$. We know that
$\psi_0(\pi_0)$ holds. By induction on $1\leq i\leq n$ select
$\pi_i$ in such a way that $\pi_{i+1}\supseteq \pi_i$, if we let
$\pi_i(\dd_i)=\dd'_i$ then $\psi_i(\pi_i,\dd'_1,\dots,\dd'_i)$
holds and $\sup(V\cap \dd_{t(i)})<\dd'_i$ and
$\Ran(\pi_i)<\dd_{t(i)}$. This is possible as $M_1,\dots,M_n$ are
elementary submodels. Finally, $\pi_n(p)$ is a condition
interlacing with $p$ by $e$ and it forces that $f(\dd'_n)=\xi$.
Now if we take the union of them plus the edge
$\{\dd_n,\dd'_n\}$ then an argument as in Lemma 7 shows that we
get a condition which forces a contradiction. \qed

\medskip
\proclaim Theorem 3. If GCH holds there is a cardinal,
cofinality, and GCH preserving (class) notion of forcing in
which for every $n$, $e$, and regular $\ll\geq\omega$ there is a
$\ll^+$-chromatic graph on $\ll^+$ all whose finite subgraphs
are subgraphs of some elements of $\KK$.

\medskip
\P For $\ll\geq\omega$ regular let $Q_\ll$ be the product of
$Q_{n,e,\ll}$ of Theorem 2 with finite supports if $\ll=\omega$,
and complete supports otherwise. Notice that $Q_\ll$ is a
$\ll^+$-c.c.~notion of forcing of cardinal $\ll^+$. For $\ll$
singular let $Q_\ll$ be the trivial forcing. 

Our notion of forcing is the Easton-support limit of the
$Q_\ll$'s, i.e., the direct limit of $P_\aa$ where $P_{\aa+1}=
P_\aa \oplus Q_\aa$ with $Q_\aa$ defined in the ground model.
For $\aa$ limit, $p\in P_\aa$ iff $p(\bb)\in Q_\bb$ for all
$\bb<\aa$, and $\vert \Dom(p)\cap\kappa \vert<\kappa$ 
for $\kappa\leq \aa$ regular.

Given $n$, $e$, and $\ll$ as in the statement of the Theorem,
the extended model can be thought as the generic extension of
some model first with $Q_{n,e,\ll}$ then with $P_\ll$ which is
of cardinal $\ll$ so it cannot change the chromatic number of a
graph from $\ll^+$ to $\ll$. 

Assume that the cofinality of some ordinal $\aa$ collapses to a
regular $\ll$. $P$ splits as $P_\ll \oplus Q_\ll \oplus R$ where
$R$ is $\leq\ll$-closed, $\vert P_\ll\vert \leq \ll$ and $Q_\ll$
is $\ll^+$-c.c., so in fact the $\ll^+$-c.c.~$P_{\ll+1}$
changes the cofinality of $\aa$ which is impossible.  This also
implies that no cardinals are collapsed. 

If $\tau$ is regular, all subsets of $\tau$ are added by the
$\tau^+$-c.c.~$P_{\tau+1}$ of cardinal $\tau^+$ so $2^\tau$
remains $\tau^+$. If $\tau$ is singular we must bound
$\tau^{\cf(\tau)}$. The sets of size $\cf(\tau)$ are added by
$P_{\cf(\tau)+1}$ so we can bound the new value of
$\tau^{\cf(\tau)}$ by $\tau^{\cf(\tau)^+}=\tau^+$. \qed
\bigskip
\centerline{{\bf References}}
\medskip
\item {[1]} P.Erd\H os, A.Hajnal: On chromatic number of graphs
and set-systems, {\sl Acta Math. Acad. Sci. Hung. {\bf 17}}
(1966), 61--99.
\item {[2]} P.Erd\H os, A.Hajnal, S.Shelah: On some general
properties of chromatic number, in: Topics in Topology, Keszthely
(Hungary), 1972, {\sl Coll. Math. Soc. J. Bolyai {\bf 8}}, 243--255.
\item {[3]} P.Erd\H os, R.Rado: Partition relations connected
with the chromatic number of graphs, {\sl
Journal of London Math. Soc. {\bf 34}} (1959), 63--72.
\item {[4]} P.Erd\H os, R.Rado: A construction of graphs without
triangles having pre-assigned order and chromatic number, {\sl
Journal of London Math. Soc. {\bf 35}} (1960), 445--448.
\item {[5]} A.Hajnal, P.Komj\'ath: What must and what need not be 
contained in a graph of
uncountable chromatic number ?, {\sl Combinatorica \bf 4} (1984),  47--52.
\item {[6]} P.Komj\'ath: The colouring number,  {\sl Proc. London Math. Soc. 
\bf
54} (1987),   
1-14.
\item {[7]} W.Taylor: Atomic compactness and elementary
equivalence, {\sl Fund. Math. {\bf 71}} (1971), 103--112.
\item {[8]} W.Taylor: Problem 42, Comb. Structures and their
applications, {\sl Proc. of the  Calgary International 
Conference}, 1969.

\bigskip

\hbox{\vbox{\hbox{P\'eter Komj\'ath}
            \hbox{Department of Computer Science}
            \hbox{E\"otv\"os University}
            \hbox{Budapest, M\'uzeum krt.~6--8}
            \hbox{1088, Hungary}
            \hbox{e-mail: {\tt kope@cs.elte.hu}}
     }
     \hskip3cm
     \vbox{\hbox{Saharon Shelah}
           \hbox{Institute of Mathematics}
           \hbox{Hebrew University}
           \hbox{Givat Ram}
           \hbox{Jerusalem, Israel}
           \hbox{e-mail: {\tt shelah@math.huji.ac.il}}}
      }

\bye